\documentclass{amsart}

\usepackage[T1]{fontenc}
\usepackage[utf8]{inputenc}
\usepackage{lmodern}
\usepackage{amsmath}
\usepackage{amssymb}
\usepackage{amsfonts}
\usepackage{amsthm}
\usepackage{tabularx}
\usepackage{graphicx}
\usepackage{tikz}
\usepackage{mathtools}
\usepackage[breaklinks]{hyperref}
\usepackage{cleveref}
\usepackage{dirtytalk}
\usepackage[hyphenbreaks]{breakurl}

\newtheorem{theorem}{Theorem}[section]
\newtheorem{corollary}[theorem]{Corollary}

\newtheorem{proposition}[theorem]{Proposition}
\theoremstyle{definition}
\newtheorem{definition}[theorem]{Definition}
\newtheorem{remark}[theorem]{Remark}
\newtheorem{example}[theorem]{Example}
\newtheorem{construction}[theorem]{Construction}

\usetikzlibrary{shapes.geometric}
\usetikzlibrary{decorations.markings}
\tikzset{middlearrow/.style={decoration={markings,mark = at position 0.5 with {\arrow{#1}},},postaction={decorate} } }

\definecolor{codegreen}{rgb}{0,0.6,0}
\definecolor{codegray}{rgb}{0.5,0.5,0.5}
\definecolor{codepurple}{rgb}{0.58,0,0.82}
\definecolor{backcolour}{rgb}{0.95,0.95,0.92}

\DeclarePairedDelimiter\abs{\lvert}{\rvert}%
\DeclarePairedDelimiter\norm{\lVert}{\rVert}%

\makeatletter
\let\oldabs\abs
\def\abs{\@ifstar{\oldabs}{\oldabs*}}
\let\oldnorm\norm
\def\norm{\@ifstar{\oldnorm}{\oldnorm*}}
\makeatother

\newcommand{\bN}{\mathbb{N}}

\newcommand{\bZ}{\mathbb{Z}}

\DeclareMathOperator{\Aut}{Aut}

\title{Finite groups as homotopy self-equivalences of finite spaces}
\author{Juan Felipe Celis-Rojas}

\begin{document}

\begin{abstract}
    We study the realization problem of finite groups as the group of homotopy classes of self-homotopy equivalences of finite spaces.
    Let $G$ be a finite group.
    Using an infinite family of pairwise non weakly homotopic asymmetric spaces we present a new construction of a finite space whose group of homotopy classes of self-homotopy equivalences is isomorphic to $G$.
\end{abstract}

\maketitle

Finite spaces are topological spaces with finitely many points. 
One can study these spaces from different perspectives such as partially ordered sets, simplicial complexes, group actions and category theory. 
This allows us to use techniques from algebraic topology and combinatorics to study these spaces. 
Thus, understanding finite spaces will give us new tools to study classical invariants from algebraic topology.

Realization problems are widely studied in algebraic topology.
Kahn's \emph{realizability problem for abstract groups} proposed in \cite{KahnRealization} studies the group of homotopy classes of self-homotopy equivalences for simply connected spaces.
An other realization problem known as \emph{Steenrod's $G$-Moore space problem}
first appeared in a topology conference in Seattle 1963m see Problem 51 in Lashoff's list  \cite{lashof_problems_1965}. 
Costoya, Gomes and Viruel presented a generalized version of Steenrod's problem for finite spaces in \cite{RealizationViaAlexandroff}. 

We focus on the realization problem for the group of homotopy classes of self-homotopy equivalences. 
Different versions of this problem have been solved.
Barmak solved it for graphs and lattices \cite{barmak_automorphism_2020}, Costoya and Viruel for elliptic spaces \cite{costoya_every_2014} and together with Alicia Tocino and Panagiote Ligouras for regular evolution algebras \cite{costoya_regular_2021} , and Chocano for Alexandroff spaces \cite{chocano_topological_2020}.
The problem we are interested in, with finite spaces, can be reduced to a graph theory problem.
The group of homotopy classes of homotopy self-equivalences of a minimal finite space is isomorphic to the automorphism group of its Hasse diagram representation.
Most results on this problem, if not all, are based on papers by Robert Frucht \cite{frucht_herstellung_1939,frucht_graphs_1949,frucht_construction_1950}. 
This article presents a new solution to the realization problem.

\hfill \break
\textbf{Theorem \ref{thm:realization}}
\textit{Let $G$ be a finite group and $S=\{g_1, \dots, g_n\}$ a set of generators of $G$. 
Then the group of homotopy classes of homotopy self-equivalences of the finite space $X(G,S)$, which is described in \Cref{def:RealizationSpace}, is isomorphic to $G$.}
\hfill \break

In \Cref{sec:finitespaces} we recall basic definitions and properties of finite spaces and briefly explain how they are related to posets.
Next we present in \Cref{sec:homotopy} the realization problem on the group of homotopy classes of self-homotopy equivalences.
Inspired by Frucht's work \cite{frucht_herstellung_1939}, we build an infinite family of non-homotopic asymmetric finite spaces, i.e. spaces with trivial group of homotopy classes of self-homotopy equivalences.
These asymmetric finite spaces are the key-point of our solution of the realization problem.
Starting with the Cayley graph of a given finite group, inserting such asymmetric finite spaces allows us to differentiate edges of distinct generators and encode the direction of all edges.
This distinction induces a restriction on the automorphism group of the space built from the Cayley graph.
It turns out that this space is a minimal finite space and then solves the realization problem.
Finally, in \Cref{sec:examples} we give a few simple but representative examples of our construction.

\hfill \break
\indent \textbf{Aknowledgements.} The author thanks the Summer in the Lab program at EPFL, Kathryn Hess and Jérôme Scherer for the guidance, Antonio Viruel for the suggestions he made on this construction and his careful reading of a draft of this preprint, and finally Peter May for motivating studying finite spaces in his REU programs.

\section{Finite topological spaces}\label{sec:finitespaces}
In this section we introduce finite topological spaces and provide some results following Peter May's approach \cite{may_finite_2003, may_finite_2003-1, may_finite_2003-2}, Barmak's book \cite{barmak_algebraic_2011} and insight from Hatcher's book \cite{hatcher_algebraic_2005}.

A \emph{finite (topological) space} is a space with finitely many points.
At first sight this definition may seem uninteresting in homotopy theorey, however finite spaces give much to talk about.
For example, for any CW-complex with finitely many cells there is a finite space which is weakly homotopic to it, see Theorem 1.7 in \cite{may_finite_2003-1}.

Notice that there are finitely many topologies on a set of $n$ points. 
In a given topology in a finite space there are finitely many open sets. 
So arbitrary intersections of open sets in a finite space are again open.
A space satisfying this condition is called an \emph{Alexandroff space}, which is abbreviated to \emph{$A$-space}.

\begin{definition}
Let $X$ be a finite space and $x \in X$. 
If $U_x$ is the intersection of all open sets containing $x$, define $\leq$ a relation on $X$ by $x \leq y$ if and only if $U_x \subseteq U_y$.
\end{definition}

\begin{remark}
This relation is reflexive and transitive. 
It is anti-symmetric if and only if $X$ is a $T_0$-space. 
If this is the case then $(X,\leq)$ is a poset.
\end{remark}

It turns out that up to homotopy equivalence, we can always assume that $X$ is $T_0$, see Proposition 1.3.1 in \cite{barmak_algebraic_2011}.
From now on we will always work with $T_0$ representatives of homotopy classes so that $(X,\leq)$ is always a poset.

In fact, there is an equivalence of categories between $A$-spaces and posets.
It will be useful to see finite spaces as posets, specially to represent them in a Hasse diagram form.

\begin{definition}
    The \emph{Hasse diagram} representation of a poset $(X,\leq)$ is a directed graph whose vertex set is $X$ and there is a directed edge $(x,y)$ for $x \ne y$ if and only if $x \leq y$ in $X$.
    We say that the level of $x \in X$ in the Hasse diagram of $X$ is
    \begin{equation*}
        level(x) = \sup_{n \in \bN}\{n : \exists x_1 < \cdots < x_n = x\}.
    \end{equation*}
\end{definition}

\begin{remark}\label{rmk:AutoHasse}
    Observe that an automorphism of a Hasse diagram must preserve the level of all vertices.
\end{remark}

\begin{definition}
Let $X$ be a finite space. 
A point $x \in X$ is called \emph{up-beat} if there is another point $y \in X\backslash\{x\}$, $x<y$ such that $x<z$ implies $y \leq z$.
A point $x \in X$ is called \emph{down-beat} if there is another point $y \in X\backslash\{x\}$, $y<x$ such that $z<x$ implies $z \leq y$.
\end{definition}

In the Hasse diagram representation of finite spaces, up-beat and down-beat points can be easily distinguished.
These are points with only one arrow going out (respectively in) them.

\begin{definition}
Let $X$ be a finite space. 
If $X$ is $T_0$ and has no up-beat or down-beat points then $X$ is called a \emph{minimal} finite space.
\end{definition}

\begin{theorem}
Let $X$ be a minimal finite space and $f: X \to X$ be a continuous map. 
If $f$ is homotopic to the identity then $f$ is the identity.
\end{theorem}
\begin{proof}
    A nice proof of this statement can be found in \cite{may_finite_2003} Theorem 6.8.
\end{proof}

\begin{corollary}\label{cor:homeo}
Let $X,Y$ be minimal finite spaces and $f: X \to Y$ be a continuous map. 
If $f$ is a homotopy equivalence then it is a homeomorphism.
\hfill $\square$
\end{corollary}

\section{Homotopy classes of self-homotopy equivalences}\label{sec:homotopy}

We will study the group of homotopy classes of self-homotopy equivalences of finite spaces, which is a homotopy invariant notion. 
Since all finite spaces are homotopic to a minimal finite space, it suffices to study it on minimal finite spaces. 
Additionally, from \Cref{cor:homeo} it follows that self-homotopy equivalences on minimal finite spaces are homeomorphisms. 
Thus the group of homotopy classes of self-homotopy equivalences of a minimal finite space is equal to its group of homeomorphisms.

\begin{remark}\label{rmk:graphiso}
    Recall from \Cref{sec:finitespaces} that a minimal finite space is also a partially ordered set. 
    So a homeomorphism of a minimal finite space corresponds to a directed graph automorphism of the Hasse-diagram of its associated poset.
\end{remark}

Our goal is to revisit a known realization problem. 
For any finite group $G$, we want to find a finite space $X$ such that its group of homotopy classes of self-homotopy equivalences is isomorphic to $G$.

\begin{remark}
    We have said before that for any finite CW-complex there is a finite space which is \emph{weakly} homotopic to it.
    It is important to keep in mind that the groups of homotopy classes of self-homotopy equivalences of these spaces are not necessarily isomorphic.
\end{remark}

\subsection{Finite spaces with trivial automorphism group}

To approach the realization problem, we choose a strategy for which we need to find an infinite family of non-isomorphic posets such that their automorphism group is trivial. 
This is inspired in the work of Frucht \cite{frucht_herstellung_1939} where he finds an infinite family of asymmetric graphs.

\begin{remark}
Notice that an isomorphism on a minimal graph induces a bijection on each level of the underlying poset.
\end{remark}

\begin{construction}\label{def:asym}
For $k \in \bN$ let $F_k$ be the finite space on $2k+8$ points such that:
\begin{enumerate}
    \item Its Hasse diagram has two levels, each with $k+4$ vertices
    \item On each level there are two vertices of degree 2 and for all $i \in \{3, \dots, k+4\}$ one vertex of degree $i$.
    \item One of the two vertices of degree 2 is connected to one vertex of degree $2$ on the other level and to the vertex of degree $k+4$ also in the other level.
    \item The other vertex of degree 2 is connected to the vertices of degree $k+3$ and $k+4$ on the other level.
\end{enumerate}
This defines $F_k$ up to isomorphism.
By construction, these finite spaces are minimal since in their Hasse diagram there are no up-beat nor down-beat points.
\end{construction}

\begin{proposition}
There are infinitely many minimal finite spaces with trivial automorphism group. In fact, the family of finite spaces $F_k$ defined just above is a family of non-homotopic asymmetric finite spaces.
\end{proposition}
\begin{proof}
We will prove that $F_k$ is asymmetric for all $k \in \bN$.

Consider the space $F_k$ and let $n=k+4$, so $F_k$ has $2n$ points. 
Recall that automorphisms can not change the degree of a vertex nor its level in the Hasse diagram. 
Thus all vertices of degree greater than 2 are fixed by all automorphisms. 
Now consider the vertices of degree two. 
One of them is connected to vertices of degree 2 and $n=k+4$, whereas the other one is connected to vertices of degree $n$ and $n-1$. 
So they can not be exchanged by any automorphism.

It follows that the only automorphism of this finite space is the identity.
\end{proof}
\begin{figure}[ht]
    \centering
    \begin{tikzpicture}
        \draw[fill] (-4,1) circle [radius = 0.075];
        \draw[fill] (-2,1) circle [radius = 0.075];
        \draw[fill] (0,1) circle [radius = 0.075];
        \draw[fill] (2,1) circle [radius = 0.075];
        \draw[fill] (4,1) circle [radius = 0.075];
        \draw[fill] (6,1) circle [radius = 0.075];
        \draw[fill] (8,1) circle [radius = 0.075];
        
        \draw[fill] (-4,-1) circle [radius = 0.075];
        \draw[fill] (-2,-1) circle [radius = 0.075];
        \draw[fill] (0,-1) circle [radius = 0.075];
        \draw[fill] (2,-1) circle [radius = 0.075];
        \draw[fill] (4,-1) circle [radius = 0.075];
        \draw[fill] (6,-1) circle [radius = 0.075];
        \draw[fill] (8,-1) circle [radius = 0.075];
        
        \draw (-4,1)--(-4,-1);
        \draw (-4,1)--(-2,-1);
        \draw (-4,-1)--(-2,1);
        
        \draw (8,1)--(-2,-1);
        \draw (8,1)--(0,-1);
        \draw (8,-1)--(-2,1);
        \draw (8,-1)--(0,1);
        
        \draw (-2,1)--(-2,-1);
        \draw (-2,1)--(0,-1);
        \draw (-2,-1)--(0,1);
        \draw (-2,1)--(2,-1);
        \draw (-2,-1)--(2,1);
        \draw (-2,1)--(4,-1);
        \draw (-2,-1)--(4,1);
        \draw (-2,1)--(6,-1);
        \draw (-2,-1)--(6,1);
        
        \draw (0,1)--(0,-1);
        \draw (0,1)--(2,-1);
        \draw (0,-1)--(2,1);
        \draw (0,1)--(4,-1);
        \draw (0,-1)--(4,1);
        \draw (0,1)--(6,-1);
        \draw (0,-1)--(6,1);
        
        \draw (2,1)--(2,-1);
        \draw (2,1)--(4,-1);
        \draw (2,-1)--(4,1);
        \draw (2,1)--(6,-1);
        \draw (2,-1)--(6,1);
        
        \draw (4,1)--(4,-1);
    \end{tikzpicture}
    \caption{The asymmetric minimal finite space $F_3$ with 14 points}
    \label{fig:asym}
\end{figure}
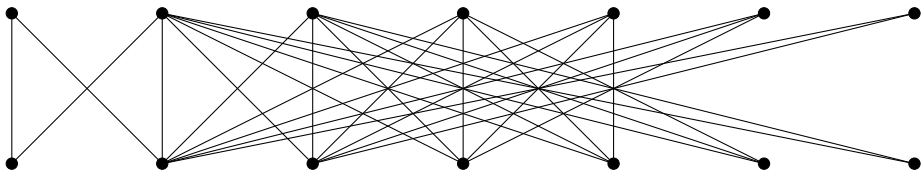

\subsection{Finite spaces with prescribed finite automorphism group}

Now we can solve the realization problem.
From now on, let $G$ be a finite group. 
We want to build a minimal finite space such that its group of automorphisms is isomorphic to $G$.

Let $S=\{g_1, \dots, g_n\}$ be a set of generators of $G$ with minimal cardinality. 
Consider $C_{G,S}$ the Cayley graph of $G$ and color the edges $(g, g_k g)$ with color $k$.
Denote this colored and directed graph by $C_{G,S}^c$.
Observe that $C_{G,S}$ has more automorphisms as a directed graph than $C_{G,S}^c$ as a colored directed graph.

\begin{construction}
Let $X$ and $B$ be $T_0$ finite spaces and $x \in X$.
A \emph{block replacement} of $x$ by $B$ on $X$ is a finite space $\widetilde{X}$ constructed in the following way from the Hasse diagram of $X$:
\begin{enumerate}
    \item Remove the vertex $x$ in the Hasse diagram of $X$ while keeping track of the edges adjacent to it;
    \item insert $B$ in the place of $x$ in the Hasse diagram of $X$; and
    \item if $e_i = (y,x)$ was an in-edge of $x$ then we add edges from $y$ to all vertices in the first level of $B$.
    Similarly, if $e_o = (x,z)$ was an out-edge of $x$ then we add edges from every vertex on the last level of $B$ to $z$.
\end{enumerate}
In this construction we call $B$ a \emph{block}, and $\widetilde{X}$ is the block replacement of $X$ at $x$ by $B$.

We can generalize this construction to build a finite space from blocks.
That is we start with a finite number of blocks and we specify how we connect them.
Whenever we connect two blocks, it means that there is an edge between all vertices on the last level of the starting block to all vertices on the first level of the ending block.
We call this \emph{construction by blocks}.
\end{construction}

\begin{construction}\label{def:RealizationSpace}
Let $G$ be a finite group with $S=\{g_1, \dots, g_n\}$ a minimal set of generators, and $C_{G,S}^c$ its Cayley graph.
We build a finite space $X(G,S)$ by block construction from the empty space adding asymmetric blocks coming from vertices and edges in $C_{G,S}^c$.

\begin{enumerate}
    \item For every element of $G$, add a block $F_0$ on the first level of $X(G,S)$;
    \item For all $1 \leq k \leq n$, for every edge of $C_{G,S}^c$ of the form $(g, g_k g)$, i.e. an edge of color $k$:
    \begin{enumerate}
        \item Add a block $F_k$ on the first level of $X(G,S)$ representing the directed edge $(g, g_k g)$;
        \item Add a block $F_{n+k}$ on the second level of $X(G,S)$, representing the starting vertex of the edge $(g, g_k g)$, and connect it the blocks of the edge $(g, g_k g)$ and vertex $g$, both in the first level of $X(G,S)$;
        \item Add a block $F_{2n+k}$ on the second level of $X(G,S)$, representing the end point of the edge $(g, g_k g)$, and connect it to the blocks of the edge $(g, g_k g)$ and vertex $g_k g$, both in the first level of $X(G,S)$.
    \end{enumerate}
\end{enumerate}
\end{construction}

\begin{remark}
The finite space $X(G,S)$ is minimal because it does not have any up-beat nor down-beat points on its Hasse diagram since every vertex has at least two adjacent edges.
\end{remark}

Understanding the Cayley graph of a group is fundamental for our solution of the realization problem.
The following theorem can be proven with Frucht's results \cite{frucht_herstellung_1939}, and has been an inspiration for many proofs on this subject.
For example, Babai studied non-colored automorphisms of the Cayley graph in \cite{BabaiAuto}.
\begin{theorem}\label{thm:colordirectedauto}
Let $G$ be a finite group. 
Then the automorphism group of the colored directed graph $C_{G,S}^c$ is isomorphic to $G$.
\end{theorem}
\begin{proof}[Sketch of proof](For a complete proof see Theorem 4-8 in \cite{ProofCayley} Chapter 4).
As $C_{G,S}^c$ is made using the Cayley graph of $G$ it is not hard to see that there is an inclusion
\begin{equation*}
    G \hookrightarrow \Aut(C_{G,S}^c)
\end{equation*}
To finish, notice that the coloring avoids having different automorphisms. 
Every automorphism is then generated by the generators of $G$. 
\end{proof}

\begin{remark}
    For a given finite group $G$ there are many  graphs whose automorphism group are isomorphic to $G$.
    Actually, Frucht proved in \cite{frucht_herstellung_1939} that there are infinitely many pairwise non-isomorphic simple graphs with automorphism group $G$.
    It is important to distinguish Frucht's graphs which are simple (undirected and uncolored), from the Cayley graph of a group which is directed and colored.
\end{remark}

\begin{theorem}\label{thm:realization}
Let $G$ be a finite group and $S=\{g_1, \dots, g_n\}$ a set of generators of $G$. 
Then the group of homotopy classes of homotopy self-equivalences of the finite space $X(G,S)$, which is described in \Cref{def:RealizationSpace}, is isomorphic to $G$.
\end{theorem}
\begin{proof}
Recall that \Cref{cor:homeo} implies that the group of homotopy classes of homotopy self-equivalences of $X(G,S)$ is isomorphic to the group of directed graph isomorphisms of its Hasse diagram as we said in \Cref{rmk:graphiso}.
Then we want to prove the following statement:
\emph{the group of directed graph isomorphisms of the Hasse diagram of $X(G,S)$ is isomorphic to $G$}.

By construction of $X(G,S)$, each automorphism of $C_{G,S}^c$ induces an automorphism on $X(G,S)$. 
Then there are inclusions
\begin{equation*}
    \Aut(C_{G,S}^c) \hookrightarrow \Aut(X(G,S)).
\end{equation*}
It remains to show the other inclusion.
Recall that an automorphism of a Hasse diagram must preserve the level of all vertices \Cref{rmk:AutoHasse}.
Let us focus only on the vertices of level 1 and 2 in $X(G,S)$.
Observing the connected components of levels 1 and 2 one detects all blocks.
The isomorphism type of each block determines whether it represents a vertex (i.e. $F_0$ isomorphism class) or an edge of color $k$ (i.e. $F_k$ isomorphism class).
An isomorphism must preserve these isomorphism classes, so it sends vertices to vertices and edges of color $k$ to edges of color $k$.

Now focus on vertices of level 3 and 4 in $X(G,S)$, the connected components are the blocks representing starting-points of edges and end-points of edges.
As above, the isomorphism type determines whether the block describes an starting/end-point of an edge and the color of the edge.
Therefore, an isomorphism not only preserves the color of edges but it also preserves their orientation.

Then every automorphism of $X(G,S)$ induces an automorphism on $C_{G,S}^c$. So there is an inclusion
\begin{equation*}
    \Aut(X(G,S)) \hookrightarrow \Aut(C_{G,S}^c).
\end{equation*}
From \Cref{thm:colordirectedauto} we know that $\Aut(C_{G,S}^c) \cong G$.
This finishes the proof by finiteness of the groups.
\end{proof}

\section{Examples}\label{sec:examples}
To illustrate the proposed solution of the realization problem, we give some simple examples of our construction.

\begin{example}
Consider the cyclic group $\bZ/3\bZ$ of order 3. 
The simplest presentation of this group is 
\begin{equation*}
    \bZ/3\bZ = \langle x \vert x^3 \rangle
\end{equation*}
Based on the presentation we can build $C_{\bZ/3\bZ,\{x\}}$ and $X(\bZ/3\bZ,\{x\})$. 
As $\bZ/3\bZ$ has only one generator there is only one color for the edges. 
We illustrate $C_{\bZ/3\bZ,\{x\}}$ in \cref{fig:CZ3}.
Notice that it is not possible to do a transposition on two vertices because it changes the direction of the edges.
\begin{figure}[ht]
    \centering
    \begin{tikzpicture}[scale=0.8]
        \draw[fill] (-2,0) circle [radius = 0.075];
        \draw[fill] (2,0) circle [radius = 0.075];
        \draw[fill] (0,3.5) circle [radius = 0.075];
        \draw[middlearrow={>}] (2,0) to [out=105, in = 315] (0,3.5);
        \draw[middlearrow={>}] (0,3.5) to [out=225, in = 75] (-2,0);
        \draw[middlearrow={>}] (-2,0) to [out=345, in = 195] (2,0);
    \end{tikzpicture}
    \caption{$C_{\bZ/3\bZ,\{x\}}$}
    \label{fig:CZ3}
\end{figure}

Next we build the finite space $X(\bZ/3\bZ,\{x\})$ as described above. 
For simplicity of the diagrams we represent asymmetric posets by vertices labeled with their type. 
We can see $X(\bZ/3\bZ,\{x\})$ in \cref{fig:XZ3}.
\begin{figure}[ht]
    \centering
    \begin{tikzpicture}
        \draw[fill] (-3,0) circle [radius = 0.075];
        \node[below] at (-3,0) {$F_0$};
        \draw[fill] (-2,0) circle [radius = 0.075];
        \node[below] at (-2,0) {$F_0$};
        \draw[fill] (-1,0) circle [radius = 0.075];
        \node[below] at (-1,0) {$F_0$};
        \draw[fill] (1,0) circle [radius = 0.075];
        \node[below] at (1,0) {$F_1$};
        \draw[fill] (2,0) circle [radius = 0.075];
        \node[below] at (2,0) {$F_1$};
        \draw[fill] (3,0) circle [radius = 0.075];
        \node[below] at (3,0) {$F_1$};
        
        \draw[fill] (-3.5,2) circle [radius = 0.075];
        \node[above] at (-3.5,2) {$F_3$};
        \draw[fill] (-2.5,2) circle [radius = 0.075];
        \node[above] at (-2.5,2) {$F_4$};
        \draw[fill] (-0.5,2) circle [radius = 0.075];
        \node[above] at (-0.5,2) {$F_3$};
        \draw[fill] (0.5,2) circle [radius = 0.075];
        \node[above] at (0.5,2) {$F_4$};
        \draw[fill] (2.5,2) circle [radius = 0.075];
        \node[above] at (2.5,2) {$F_3$};
        \draw[fill] (3.5,2) circle [radius = 0.075];
        \node[above] at (3.5,2) {$F_4$};
        
        \draw (1,0) -- (-3.5,2);
        \draw (1,0) -- (-2.5,2);
        \draw (2,0) -- (-0.5,2);
        \draw (2,0) -- (0.5,2);
        \draw (3,0) -- (2.5,2);
        \draw (3,0) -- (3.5,2);
        
        \draw (-3,0) -- (-3.5,2);
        \draw (-3,0) -- (3.5,2);
        \draw (-2,0) -- (-2.5,2);
        \draw (-2,0) -- (-0.5,2);
        \draw (-1,0) -- (0.5,2);
        \draw (-1,0) -- (2.5,2);
    \end{tikzpicture}
    \caption{$X(\bZ/3\bZ,\{x\})$}
    \label{fig:XZ3}
\end{figure}
\end{example}

\begin{example}
Consider the dihedral group $D_6$ of order 6. One presentation of this group is:
\begin{equation*}
    D_6 = \langle \tau, \sigma \vert \tau^3, \sigma^2, \tau\sigma\tau\sigma^{-1} \rangle
\end{equation*}
Based on this presentation we can draw the colored Cayley graph of $D_6$, shown in \cref{fig:CD6}. 
The graph $C_{D_6,\{\tau,\sigma\}}$ has two colors, one for each generator.
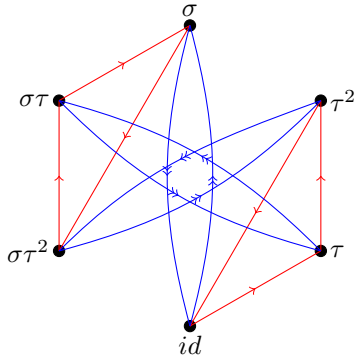
\begin{figure}[ht]
    \centering
    \begin{tikzpicture}
        \draw[fill] (0,-2) circle [radius=0.075];
        \node [below] at (0,-2) {$id$};
        \draw[fill] (1.74,-1) circle [radius=0.075];
        \node [right] at (1.74,-1) {$\tau$};
        \draw[fill] (1.74,1) circle [radius=0.075];
        \node [right] at (1.74,1) {$\tau^2$};
        \draw[fill] (0,2) circle [radius=0.075];
        \node [above] at (0,2) {$\sigma$};
        \draw[fill] (-1.74,1) circle [radius=0.075];
        \node [left] at (-1.74,1) {$\sigma\tau$};
        \draw[fill] (-1.74,-1) circle [radius=0.075];
        \node [left] at (-1.74,-1) {$\sigma\tau^2$};
        
        \draw[middlearrow={>},red] (0,-2)--(1.74,-1);
        \draw[middlearrow={>},red] (1.74,-1)--(1.74,1);
        \draw[middlearrow={>},red] (1.74,1)--(0,-2);
        
        \draw[middlearrow={>},red] (0,2)--(-1.74,-1);
        \draw[middlearrow={>},red] (-1.74,-1)--(-1.74,1);
        \draw[middlearrow={>},red] (-1.74,1)--(0,2);
        
        \draw[middlearrow={>>},blue] (0,-2) to [out=75, in = 285] (0,2);
        \draw[middlearrow={>>},blue] (0,2) to [out=255, in = 105] (0,-2);
        
        \draw[middlearrow={>>},blue] (1.74,-1) to [out=135, in = 345] (-1.74,1);
        \draw[middlearrow={>>},blue] (-1.74,1) to [out=315, in = 165] (1.74,-1);
        
        \draw[middlearrow={>>},blue] (1.74,1) to [out=200, in = 45] (-1.74,-1);
        \draw[middlearrow={>>},blue] (-1.74,-1) to [out=15, in = 225] (1.74,1);
    \end{tikzpicture}
    \caption{$C_{D_6,\{\tau,\sigma\}}$}
    \label{fig:CD6}
\end{figure}

Next we build $X(D_6,\{\tau,\sigma\})$. 
The process we described above gives us the finite space in \cref{fig:XD6}.
\begin{figure}[ht]
    \centering
    \begin{tikzpicture}
        \draw[fill] (-5,0) circle [radius = 0.075];
        \node[below] at (-5,0) {\tiny{$F_0$}};
        \draw[fill] (-4.5,0) circle [radius = 0.075];
        \node[below] at (-4.5,0) {\tiny{$F_0$}};
        \draw[fill] (-4,0) circle [radius = 0.075];
        \node[below] at (-4,0) {\tiny{$F_0$}};
        \draw[fill] (-3.5,0) circle [radius = 0.075];
        \node[below] at (-3.5,0) {\tiny{$F_0$}};
        \draw[fill] (-3,0) circle [radius = 0.075];
        \node[below] at (-3,0) {\tiny{$F_0$}};
        \draw[fill] (-2.5,0) circle [radius = 0.075];
        \node[below] at (-2.5,0) {\tiny{$F_0$}};
        
        \draw[fill] (-1.5,0) circle [radius = 0.075];
        \node[below] at (-1.5,0) {\tiny{$F_1$}};
        \draw[fill] (-1,0) circle [radius = 0.075];
        \node[below] at (-1,0) {\tiny{$F_1$}};
        \draw[fill] (-0.5,0) circle [radius = 0.075];
        \node[below] at (-0.5,0) {\tiny{$F_1$}};
        \draw[fill] (0,0) circle [radius = 0.075];
        \node[below] at (0,0) {\tiny{$F_1$}};
        \draw[fill] (0.5,0) circle [radius = 0.075];
        \node[below] at (0.5,0) {\tiny{$F_1$}};
        \draw[fill] (1,0) circle [radius = 0.075];
        \node[below] at (1,0) {\tiny{$F_1$}};
        
        \draw[fill] (2,0) circle [radius = 0.075];
        \node[below] at (2,0) {\tiny{$F_2$}};
        \draw[fill] (2.5,0) circle [radius = 0.075];
        \node[below] at (2.5,0) {\tiny{$F_2$}};
        \draw[fill] (3,0) circle [radius = 0.075];
        \node[below] at (3,0) {\tiny{$F_2$}};
        \draw[fill] (3.5,0) circle [radius = 0.075];
        \node[below] at (3.5,0) {\tiny{$F_2$}};
        \draw[fill] (4,0) circle [radius = 0.075];
        \node[below] at (4,0) {\tiny{$F_2$}};
        \draw[fill] (4.5,0) circle [radius = 0.075];
        \node[below] at (4.5,0) {\tiny{$F_2$}};
        
        \draw[fill] (-6,2) circle [radius = 0.075];
        \node[above] at (-6,2) {\tiny{$F_3$}};
        \draw[fill] (-5.5,2) circle [radius = 0.075];
        \node[above] at (-5.5,2) {\tiny{$F_4$}};
        \draw[fill] (-5,2) circle [radius = 0.075];
        \node[above] at (-5,2) {\tiny{$F_3$}};
        \draw[fill] (-4.5,2) circle [radius = 0.075];
        \node[above] at (-4.5,2) {\tiny{$F_4$}};
        \draw[fill] (-4,2) circle [radius = 0.075];
        \node[above] at (-4,2) {\tiny{$F_3$}};
        \draw[fill] (-3.5,2) circle [radius = 0.075];
        \node[above] at (-3.5,2) {\tiny{$F_4$}};
        \draw[fill] (-3,2) circle [radius = 0.075];
        \node[above] at (-3,2) {\tiny{$F_3$}};
        \draw[fill] (-2.5,2) circle [radius = 0.075];
        \node[above] at (-2.5,2) {\tiny{$F_4$}};
        \draw[fill] (-2,2) circle [radius = 0.075];
        \node[above] at (-2,2) {\tiny{$F_3$}};
        \draw[fill] (-1.5,2) circle [radius = 0.075];
        \node[above] at (-1.5,2) {\tiny{$F_4$}};
        \draw[fill] (-1,2) circle [radius = 0.075];
        \node[above] at (-1,2) {\tiny{$F_3$}};
        \draw[fill] (-0.5,2) circle [radius = 0.075];
        \node[above] at (-0.5,2) {\tiny{$F_4$}};
        
        \draw[fill] (0.5,2) circle [radius = 0.075];
        \node[above] at (0.5,2) {\tiny{$F_5$}};
        \draw[fill] (1,2) circle [radius = 0.075];
        \node[above] at (1,2) {\tiny{$F_6$}};
        \draw[fill] (1.5,2) circle [radius = 0.075];
        \node[above] at (1.5,2) {\tiny{$F_5$}};
        \draw[fill] (2,2) circle [radius = 0.075];
        \node[above] at (2,2) {\tiny{$F_6$}};
        \draw[fill] (2.5,2) circle [radius = 0.075];
        \node[above] at (2.5,2) {\tiny{$F_5$}};
        \draw[fill] (3,2) circle [radius = 0.075];
        \node[above] at (3,2) {\tiny{$F_6$}};
        \draw[fill] (3.5,2) circle [radius = 0.075];
        \node[above] at (3.5,2) {\tiny{$F_5$}};
        \draw[fill] (4,2) circle [radius = 0.075];
        \node[above] at (4,2) {\tiny{$F_6$}};
        \draw[fill] (4.5,2) circle [radius = 0.075];
        \node[above] at (4.5,2) {\tiny{$F_5$}};
        \draw[fill] (5,2) circle [radius = 0.075];
        \node[above] at (5,2) {\tiny{$F_6$}};
        \draw[fill] (5.5,2) circle [radius = 0.075];
        \node[above] at (5.5,2) {\tiny{$F_5$}};
        \draw[fill] (6,2) circle [radius = 0.075];
        \node[above] at (6,2) {\tiny{$F_6$}};
        
        \draw (-1.5,0)--(-6,2);
        \draw (-1.5,0)--(-5.5,2);
        
        \draw (-1,0)--(-5,2);
        \draw (-1,0)--(-4.5,2);
        
        \draw (-0.5,0)--(-4,2);
        \draw (-0.5,0)--(-3.5,2);
        
        \draw (0,0)--(-3,2);
        \draw (0,0)--(-2.5,2);
        
        \draw (0.5,0)--(-2,2);
        \draw (0.5,0)--(-1.5,2);
        
        \draw (1,0)--(-1,2);
        \draw (1,0)--(-0.5,2);
        
        \draw (2,0)--(0.5,2);
        \draw (2,0)--(1,2);
        
        \draw (2.5,0)--(1.5,2);
        \draw (2.5,0)--(2,2);
        
        \draw (3,0)--(2.5,2);
        \draw (3,0)--(3,2);
        
        \draw (3.5,0)--(3.5,2);
        \draw (3.5,0)--(4,2);
        
        \draw (4,0)--(4.5,2);
        \draw (4,0)--(5,2);
        
        \draw (4.5,0)--(5.5,2);
        \draw (4.5,0)--(6,2);
        
        \draw (-5,0)--(-6,2);
        \draw (-5,0)--(-3.5,2);
        \draw (-5,0)--(0.5,2);
        \draw (-5,0)--(2,2);
        
        \draw (-4.5,0)--(-5.5,2);
        \draw (-4.5,0)--(-5,2);
        \draw (-4.5,0)--(2.5,2);
        \draw (-4.5,0)--(4,2);
        
        \draw (-4,0)--(-4.5,2);
        \draw (-4,0)--(-4,2);
        \draw (-4,0)--(4.5,2);
        \draw (-4,0)--(6,2);
        
        \draw (-3.5,0)--(-3,2);
        \draw (-3.5,0)--(-0.5,2);
        \draw (-3.5,0)--(1,2);
        \draw (-3.5,0)--(1.5,2);
        
        \draw (-3,0)--(-2.5,2);
        \draw (-3,0)--(-2,2);
        \draw (-3,0)--(3,2);
        \draw (-3,0)--(3.5,2);
        
        \draw (-2.5,0)--(-1.5,2);
        \draw (-2.5,0)--(-1,2);
        \draw (-2.5,0)--(5,2);
        \draw (-2.5,0)--(5.5,2);
    \end{tikzpicture}
    \caption{$X(D_6,\{\tau,\sigma\})$}
    \label{fig:XD6}
\end{figure}
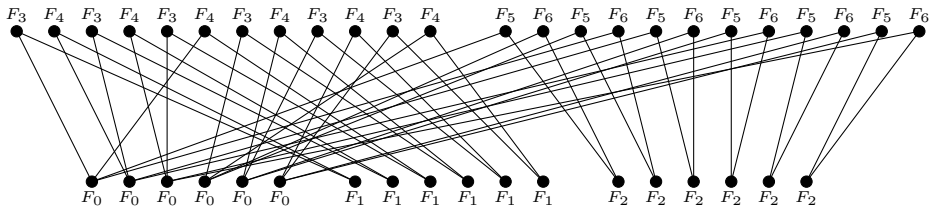
\end{example}

\bibliographystyle{alpha}

\end{document}